\documentclass{amsart}
\usepackage[utf8]{inputenc}
\usepackage{amsmath}
\usepackage{amssymb}
\usepackage{amsfonts}
\usepackage{cleveref}
\usepackage{tikz-cd}
\usepackage{stmaryrd}

\newtheorem{theorem}{Theorem}[section]

\theoremstyle{definition}

\newcommand{\ilimit}{\mbox{$\,\displaystyle{\lim_{\longleftarrow}}\,$}}
\usepackage[arrow,matrix,curve]{xy}\SilentMatrices
\def\xyma{\xymatrix@M.7em}

\begin{document}

\title{Homotopy patterns in group theory}
\author{Roman Mikhailov}

\address{Laboratory of Modern Algebra and Applications, St. Petersburg State University, 14th Line, 29b,
Saint Petersburg, 199178 Russia and St. Petersburg Department of Steklov Mathematical Institute, Fontanka 27, Saint Petersburg, 191023 Russia} \email{romanvm@mi.ras.ru}

\begin{abstract}
This is a survey. The main subject of this survey is the homotopical or homological nature of certain structures which appear in classical problems about groups, Lie rings and group rings. It is well known that the (generalized) dimension subgroups have complicated combinatorial theories. In this paper we show that, in certain cases, the complexity of these theories is based on homotopy theory. The derived functors of non-additive functors, homotopy groups of spheres, group homology etc appear naturally in problems formulated in purely group-theoretical terms. The variety of structures appearing in the considered context is very rich. In order to illustrate it, we present this survey as a trip passing through examples having a similar nature.  
\end{abstract}

\thanks{
The work is supported by the grant of the Government of the Russian Federation for the state
support of scientific research carried out under the supervision of leading scientists, agreement
14.W03.31.0030 dated 15.02.2018.1.}

\maketitle

\section{Introduction} It would be nice to discuss the title first. What will we mean by "group theory"? Obviously not a collection of funny stories on various constructions of groups with exotic properties. As well as not a classification of groups from the point of view of certain structure theory. In this paper, the groups will be the main category we will work in. In some cases it will be changed by Lie rings over integers. So, for us the "group theory" will mean just a (mostly functorial) life inside the category of groups. 

We will not define the \textit{homotopy pattern} as a concept. The term \textit{pattern} itself is multi-valued. It is used in different philosophical contexts and usually is understood intuitively. 
One can try to define it after reading this paper. Somebody can understand it as a \textit{system of signs} of homotopy origin or a \textit{collection of relations} that comes from the homotopy theory. 

The main claim or "formula" which we state here is the following abstract relation:
$$
(\Phi):\ \ \ \ \ \ \ \ \ \ \ \ \ \ \ \ \  \frac{\textsf{Intersection\ of\ subsctructures}}{\textsf{Obvious\ part}}\supset {\textsf{Homotopy\ pattern}},
$$
We will take some structure like a group, group ring, Lie ring or universal enveloping algebra, consider certain substructures, their intersection, take its quotient by some obvious part and see that, in many cases, this quotient contains elements of homotopical or homological nature. An \textit{obvious part} is not defined in an unified way, usually it is a maximal substructure of the intersection defined using given type of operations (for example, it is not an intersection itself). As a rule, the \textit{obvious part} is a more explicit construction than the \textit{intersection of substructures}. One can consider the left hand site of the formula $(\Phi)$ as "implicity modulo explicity", or define a hierarchy of the explicity and take the quotients of its terms.  

It is a time to stop saying general words and do some math. 
Let $G$ be a group, $\mathbb Z[G]$ its integral group ring, $\textbf b$ a two-sided ideal
of $\mathbb Z[G]$. The problem of identification of the subgroup 
$$
D(G, \textbf{b}) :=
G \cap (1 + \textbf{b}) = \langle g\in G\ |\ g-1 \in \textbf{b}\rangle$$ 
is a fundamental problem in the theory
of groups and group rings. It is often the case that a certain normal subgroup
$N(G, \textbf{b})$ of $G$, is easily seen to be contained in $D(G, \textbf{b})$ and explicitly defined in terms of $G$ only, without using the group ring (and is the largest subgroup with such property). Computation of the quotient $D(G, \textbf{b})/N(G, \textbf{b})$ usually becomes a challenging problem. This case will be the first example of our formula $(\Phi)$. Various choices of the ideal $\textbf{b}$ lead to the derived functors inside the quotients $D(G, \textbf{b})/N(G, \textbf{b})$. We will discuss the derived functors and their appearance in this context in Section 2. 

The main examples of the above type are classical dimension subgroups. Let $\textbf{g}$ be the augmentation ideal of
$\mathbb Z[G]$. The subgroups $D_n(G):= G \cap 
(1 + \textbf{g}^n),\ n\geq 1$ are known as
\textit{dimension subgroups}. It is easy to see that, for any $G$ and $n\geq 1$, the
dimension subgroups contain the terms of lower central series of $G$:
$\gamma_n(G)\subseteq D_n(G)$. The lower central series are defined inductively as follows $\gamma_1(G):=G,\ \gamma_{n+1}(G):=[\gamma_n(G),G]=\langle [x,y],\ x\in \gamma_n(G), y\in G\rangle^G,\ n\geq 1$. Is it true that $D_n(G)=\gamma_n(G)$? This problem was open during many years (see Section 5 for some history). Our formula $(\Phi)$ states that the dimension quotients $D_n/\gamma_n$ may contain certain elements of homotopical nature. This is exactly the case. Section 5 is about it. In particular, we will show in details in Section 5, how an element $\mathbb Z/3\subset \pi_6(S^2)$ is related to the 3-torsion in dimension quotients for Lie rings. 

The formula $(\Phi)$, which presents not a rigorous statement but feeling, goes through the paper. In Section 3 we will see that the homotopy groups of certain spaces can be described as intersections of subgroups in a group modulo a natural commutator subgroup. In particular, the homotopy groups of the 2-sphere are given in this way. This is the well-known Wu formula and its variations. Sections 4 is about combinatorics of Wu-type formulas. One can see examples of \textit{homotopy patterns} as combinations of brackets in groups and Lie rings. As we mentioned already, Section 5 shows how these combinations can be applied to the classical dimension problem.  

Section 6 is a bit isolated, however, it's subject fits smoothly in a general context of the paper. Section 6 is about the method of construction of functors via (derived) limits. We show how to present certain derived functors, group homology, the forth dimension quotient via limits over the category of group presentations. So the theory of limits becomes a unified theory for different functors discussed in this paper. At the end of Section 6 we briefly discuss so-called $\textbf{fr}$-language, a combinatorial-linguistic game which can be used in the study of functors.  

On October 2, 2021, my friend and teacher Inder Bir Singh Passi passed away. For all those 19 years that we were in contact, I was deeply touched by his delicacy, kindness, and empathy for other people. In 2002, he invited me to visit India for the first time. That visit changed my life. I dedicate this text to his memory.
\section{Derived functors}
The derived
functors in the sense of Dold-Puppe \cite{DP:1961} are defined as follows. For an
abelian group $A$ and an endofunctor $F$ on the category of
abelian groups, the derived functor of $F$ is given as
$$
L_iF(A,n)=\pi_i(FKP_\ast[n]),\ i\geq 0, n\geq 0
$$
where $P_\ast \to A$ is a projective resolution of $A$, and
 $K$ is  the Dold-Kan transform,  inverse to the Moore normalization
 functor from simplicial abelian groups to chain complexes. For the simplicity we write $L_iF(A):=L_iF(A,0).$ 

Derived functors appear naturally in the theory of Eilenberg-MacLane spaces, Moore spaces and general homotopy theory. For example, for an abelian group $A$, homology $H_*(A)$ can be filtred in a way that the graded pieces are derived functors of the exterior powers $L_i\Lambda^j(A)$ \cite{Breen}. In particular, there exist the following natural exact sequences 
\begin{align*}
& 0\to \Lambda^3(A)\to H_3(A)\to L_1\Lambda^2(A)\to 0\\ 
& 0\to \Lambda^4(A)\to H_4(A)\to L_1\Lambda^3(A)\to 0. 
\end{align*}
These exact sequences split as sequences of abelian groups but do not split naturally.  This is a common situation, the homology and homotopy functors usually present non-trivial gluing of derived functors of different type. 

As a rule, the derived functors have a complicated structure. Here we will describe a couple of them, in a sense the simplest ones, and show how do they appear in the context of group rings. First let us recall the well-known description (see, for example,
MacLane \cite{Mac1}), of the derived functor of the tensor square
$L_1\otimes^2(A)=\textsf{Tor}(A,A)$. Given an abelian group $A$, the group $\textsf{Tor}(A,A)$
is generated by the $n$-linear expressions $\tau_h(a_1,a_2)$
(where all $a_i$ belong to the subgroup $ {}_hA:=\{a\in A\,|\,
ha=0\},  \ h
>0)$, subject to the so-called slide relations
\begin{equation}
\label{slide} \tau_{hk}(a_1,\,a_2) = \tau_{h}(ka_1, \,a_2)
\end{equation} for all $i$ whenever  $hka_1 = 0$ and $ha_2=0$ and
analogous relation, where the roles of $a_1,\ a_2$ are changed.

The natural projection of the tensor square to the symmetric square $\otimes^2 \to {\textsf S}^2$  induces a
natural epimorphism
$\label{epi} L_1\otimes^2(A)\to L_1{\textsf S}^2(A)$ which maps
the generator $\tau_h(a_1,\,a_2)$ of $ L_1\otimes^2(A)=\textsf{Tor}(A,\,A)$ to the generator $\beta_h(a_1,\,a_2)$ of
$L_1{\textsf S}^2(A)$ so that the kernel of this map is
generated by the elements $\tau_h(a,\,a),\ a\in\ _hA.$

The functor $L_1{\textsf S}^2$ appears as a quadratic piece of the homology of Eilenberg-MacLane spaces $H_5K(-,2)$: 
$$
0\to L_1{\textsf S}^2(A)\to H_5K(A,2)\to \textsf{Tor}(A,\mathbb Z/2)\to 0. 
$$

It is shown by Jean  in  \cite{Jean} that the first derived functor of the symmetric cube can be described as follows
\begin{equation}\label{syder}
L_1{\textsf S}^{3}(A)\simeq
(L_1{\textsf S}^{2}(A)\otimes A)/Jac_{S},
\end{equation}
where $Jac_{S}$ is the subgroup generated by elements of the form (Jacobi-type elements)
$$
\beta_h(x_1,\,x_2)\otimes x_3+\beta_h(x_1,\,x_3)\otimes
x_2+\beta_h(x_2,\,x_3)\otimes x_1
$$
with $x_i \in {}_hA$.

Recall one more property of the functor $L_1{\textsf S}^2$. Suppose that an abelian group $A$ is presented as a quotient $A=Q/U$, for a free abelian $Q$ and its subgroup $U$. Then ${\textsf S}^2(A), L_1{\textsf S}^2(A)$ are naturally isomorphic to the zeroth and first homology of the Kozsul-type complex \cite{Kock:2001} 
$$
\Lambda^2(U)\to U\otimes Q\to {\textsf S}^2(Q). 
$$
The maps in this sequence are natural and easily can be recognized. Now consider the following diagram with exact columns
\begin{equation}\label{firstdia}
\xyma{& \Lambda^2(U)\ar@{>->}[r]\ar@{>->}[d] & U\otimes
Q\ar@{->}[r]\ar@{>->}[d] & \mathsf
S^2(Q)\ar@{=}[d]\\
& \Lambda^2(Q)\ar@{>->}[r]\ar@{->>}[d] & Q\otimes Q\ar@{->}[r]\ar@{->>}[d] & \mathsf S^2(Q)\\
L_1{\textsf S}^2(A)\ar@{>->}[r] &
\Lambda^2(Q)/\Lambda^2(U)\ar@{->}[r] & Q/U\otimes Q\ar@{->>}[r] &
\textsf{S}^2(A) }
\end{equation}
All vertical maps in this diagram are obvious, the middle horizontal sequence is the classical Kozsul short exact sequence. The lower sequence is exact.   

Now we return to the theory of nonabelian groups. Let $F$ be a free group, $R$ its normal subgroup, ${\textbf f}=(F-1)\mathbb Z[F],$ ${\textbf r}=(R-1)\mathbb Z[F]$, as always. Let us describe the generalized dimension subgroup $F\cap (1+\textbf{rf}+{\textbf f}^3)$. Obviously this subgroup contains $\gamma_2(R)\gamma_3(F)$, however, this is not a complete description. In order to find the remaining part, denote $Q:=F_{ab},\ U:=R/R\cap [F,F]$. Observe now that there are natural vertical isomorphisms
$$
\xyma{\Lambda^2(Q)/\Lambda^2(U)\ar@{->}[r] \ar@{=}[d] & Q/U\otimes Q\ar@{=}[d]\\ \frac{\gamma_2(F)}{\gamma_2(R)\gamma_3(F)}\ar@{->}[r] & \frac{{\textbf f}^2}{\textbf{rf}+{\textbf f}^3}}
$$ 
and the lower map is induced by $g\mapsto g-1$. That is, the needed generalized dimension quotient can be described as (see \cite{HMP} for another proof of this statement)
\begin{equation}\label{l1s}
\frac{F\cap (1+\textbf{rf}+{\textbf f}^3)}{\gamma_2(R)\gamma_3(F)}=L_1{\textsf S}^2(G_{ab}), 
\end{equation}
where $G=F/R.$
It is easy to lift the element from $L_1{\textsf S}^2(F/R\gamma_2(F))$ to $F\cap (1+\textbf{rf}+{\textbf f}^3)$. Given $\beta_h(a_1,a_2)$, denote by $f_1,f_2$ the preimages of $a_1,a_2$ in $F$, then $f_1^h,f_2^h\in R\gamma_2(F)$. Now the needed image is given as $[f_1,f_2]^h$. 

The identification (\ref{l1s}) is an example of a situation described in our formula $(\Phi).$ The subgroup $\gamma_2(R)\gamma_3(F)$ is the maximal obvious subgroup of the intersection $F\cap (1+\textbf{rf}+{\textbf f}^3).$ In general, for two ideals $\textbf a, \textbf b,$ the maximal obvious part of $D(F,\textbf a+\textbf b)$ is the product $D(F,\textbf a)D(F,\textbf b).$ In the case above, $D(F,\textbf{rf})=\gamma_2(R)$ and $D(F,\textbf f^3)=\gamma_3(F).$ For the situation of arbitrary ideals ${\bf a}$ and ${\bf b}$, the quotient
$$
\frac{D(F,{\textbf a}+{\textbf b})}{D(F,{\textbf a})D(F,{\textbf b})}
$$
shows how these ideals are "linked" in $\mathbb Z[F]$ and in many cases has a homological description, what agrees with our formula $(\Phi).$

Here is one more example related to the functor $L_1{\textsf S}^2.$ 

Recall the so-called Fox subgroup problem (\cite{Fox}, page 557; \cite{Birman:1974}, Problem 13; \cite{Gupta:1987}). It asks for the identification of the normal subgroup $F(n,\,R):=F\cap (1+{\textbf r}{\textbf f}^n)$ for a  free group $F$ and its normal subgroup $R$. A solution to this problem has been given by I. A. Yunus \cite{Yunus:1984} and Narain Gupta (\cite{Gupta:1987}, Chapter III). It turns out that while $F(1,\,R)=\gamma_2(R)$, $F(2,\,R)=[R\cap \gamma_2(F),\,R\cap \gamma_2(F)]\gamma_3(R)$, the identification of $F(n,\,R), n\geq 3$, is given as an isolator of a subgroup. For instance, $F(3,\,R)=\sqrt{G(3,\,R)}$, where $$G(3,\,R):=\gamma_2(R\cap \gamma_3(F))[[R\cap \gamma_2(F),\,R],\, R\cap \gamma_2(F)]\gamma_4(R).$$ 

It is shown in \cite{MP2019} that there is a natural isomorphism
$$
\frac{F(3,\,R)}{G(3,\,R)}\simeq L_1{\textsf S}^2\left(\frac{R\cap \gamma_2(F)}{\gamma_2(R)(R\cap \gamma_3(F))}\right).
$$
 
As in the simplest example, the derived functor $L_1{\textsf S}^2$ can be lifted to the generalized dimension subgroup. This gives a complete description of the third Fox subgroup:
$F(3,\,R)=G(3,\,R)W$,
where $W$ is a subgroup of $F$, generated by elements
$$
[x,\,y]^m[x,\,s_y]^{-1}[y,\,s_x]
$$
with
\begin{align*}
& x^m=r_xs_x,\ r_x\in R\cap \gamma_3(F),\ s_x\in \gamma_2(R)\\
& y^m=r_ys_y,\ r_y\in R\cap \gamma_3(F),\ s_y\in \gamma_2(R).
\end{align*}
Higher generalizations of this description are of interest. For instance, the complete description of the fourth Fox subgroup $F\cap (1+\textbf{rf}^4)$ should be related to the derived functors of certain cubical functors.   

The above examples lie in the tip of an iceberg, and present the simplest illustration of deep relation between generalized dimension subgroups and derived functors. The derived functors of high degree polynomial functors as well as subgroups defined by ideals have complicated structure. In many cases the same kind of tricks as above like diagram chasing together with group theoretical identifications lead to surprising connections. 

Here are some other examples of descriptions of the generalized dimensions subgroups which use the derived functors [we assume that $G=F/R$]:
\begin{align*}
& \frac{F\cap (1+{\textbf r}^2{\textbf f}+{\textbf f}^4)}{\gamma_3(R)\gamma_4(F)}=L_2 {\textsf L}_s^3(G_{ab})\\ 
& \frac{F\cap (1+{\textbf f}(F'-1)+\textbf{rf}^3+{\textbf f}^4)}{[\gamma_2(R),F]\gamma_4(F)}=L_1{\textsf S}^3(G_{ab})\\ 
& \frac{F\cap (1+\textbf{frf}+{\textbf f}^4)}{[\gamma_2(R),F]\gamma_4(F)}=L_1{\textsf S}^3(G_{ab})\ \ \ \ \ \ \text{[provided}\ G_{ab}\ \text{is 2-torsion-free]}.
\end{align*}
For an abelian group $A$ the third super-Lie functor ${\textsf L}_s^3(A)$ is generated by brackets $\{a,b,c\},$ $a,b,c\in A$ which are additive on each variable with the following relations: $\{a,b,c\}=\{b,a,c\},$ $\{a,b,c\}+\{c,a,b\}+\{b,c,a\}=0$. The derived functors of higher super-Lie functors appear in the description of the subgroup $F\cap (1+{\textbf r}^n{\textbf f}+{\textbf f}^{n+2})$ for $n\geq 2$. 

Some exotic examples of the same nature can be found in \cite{MP:2016}. In particular, there are the following descriptions of the generalized dimension subgroups:
\begin{align*}
& \frac{F\cap (1+\textbf{rfr}+\textbf{sr})}{\gamma_2(S)\gamma_3(R)}=L_1\textsf S^2(H_2(G)),\\ 
& \frac{F\cap (1+\textbf s^2\textbf r+\textbf r^2\textbf{fr})}{\gamma_3(S)\gamma_4(R)}=L_2{\textsf L}_s^3(H_2(G)),
\end{align*}
where $S=[F,R],$ $\textbf s=(S-1)\mathbb Z[F]$ and $H_2(G)$ the second integral group homology. 

\section{Homotopy pushouts}
In this Section we will see that, not only derived functors, but also the homotopy groups of certain spaces can be presented in group theoretical terms. We start with the simplest case. 
Let $G$ be a group, $R,S$ its normal subgroups. Consider a homotopy pushout 
$$
\xyma{K(G,1)\ar@{->}[r] \ar@{->}[d] & K(G/R,1)\ar@{->}[d] \\ K(G/S,1) \ar@{->}[r] & X}
$$ 
where the maps between classifying spaces $K(-,1)$ are induced by natural epimorphisms $G\to G/R,\ G\to G/S$. Classical van Kampen theorem implies that $\pi_1(X)=G/RS$. In \cite{Brown}, the second homotopy group of $X$ is described as follows: 
$$
\pi_2(X)\simeq\frac{R\cap S}{[R,S]}. 
$$
This result gives a topological interpretation of the difference between the intersection and the commutator of a pair of subgroups. [For a pair $R,S$ of normal subgroups of a free group $F$, $[R,R]\cap [S,S]\subseteq [R,S]$. One can prove it as an exercise.]  

A fact that the second homotopy group is related to certain group-theorical construction is not surprising. The theory of (non)aspherical group presentations, identity sequences etc is about the properties of the second homotopy group of a standard complex constructed from a given group presentation. Next we will show how to extend the above result for the case of three normal subgroups. 

Let $R,S,T$ be normal subgroups of $G$. Define an analog of the commutator subgroup as 
$$
[\![R,S,T]\!]:=[R\cap S, T][S\cap T, R][T\cap R,S].
$$
Consider a homotopy pushout 

$$
{\small { \xyma{& & K(G,1) \ar@{->}[rr] \ar@{->}[ldd]
\ar@{-}[dd] && K(G/R,1) \ar@{->}[ddd] \ar@{->}[ldd]\\ \\
& K(G/S,1) \ar@{->}[ddd] \ar@{->}[rr] & \ar@{->}[d] & K(G/RS,1) \ar@{->}[ddd]\\
& & K(G/T,1) \ar@{-}[r] \ar@{->}[ldd] & \ar@{->}[r] & K(G/RT,1)
\ar@{->}[ldd]\\ & & \\ & K(G/ST,1) \ar@{->}[rr] & & X  } } }
$$ 
The lower homotopy groups of $X$ are described as follows (see \cite{EllisMikhailov}):
\begin{align*}
& \pi_1(X)\simeq G/RST,\\ 
& \pi_2(X)\simeq\frac{RS\cap RT}{R(S\cap T)},\\ 
& \pi_3(X)\simeq\frac{R\cap S\cap T}{[\![R,S,T]\!]}. 
\end{align*} 
At first glance it might seem that $\frac{RS\cap RT}{R(S\cap T)}$ is not symmetric on $R,S,T$ and that the subgroup $R$ plays a special role. However, the homotopy pushout $X$ is symmetric and the above description of $\pi_2$ provides a proof of the following isomorphisms 
$$
\frac{RS\cap RT}{R(S\cap T)}\simeq \frac{SR\cap ST}{S(R\cap T)}\simeq \frac{TR\cap TS}{T(R\cap S)}. 
$$  

The above description shows that the third homotopy group appears quite naturally in the context of group theory. Next we will show how to get a purely group-theoretic result using given homotopy identifications. 

For a triple of normal subgroup $R,S,T$ of $G,$ consider the corresponding ideals in $\mathbb Z[G]:$ ${\textbf r}:=(R-1)\mathbb Z[G],$ ${\textbf s}:=(S-1)\mathbb Z[G],$ ${\textbf t}:=(T-1)\mathbb Z[G].$ An obvious ring-theoretic analog of the subgroup $[\![R,S,T]\!]$ is the following ideal:
$$
(\!({\textbf r},{\textbf s},{\textbf t})\!):={\textbf r}({\textbf s}\cap {\textbf t})+({\textbf s}\cap {\textbf t}){\textbf r}+{\textbf s}({\textbf t}\cap {\textbf r})+({\textbf t}\cap {\textbf r}){\textbf s}+{\textbf t}({\textbf r}\cap {\textbf s})+({\textbf r}\cap {\textbf s}){\textbf t}.
$$
It is easy to check that, for any $w\in [\![R,S,T]\!],$ $w-1\in (\!({\textbf r},{\textbf s},{\textbf t})\!)$ and one ask about the structure of the quotient
$$
\frac{G\cap (1+(\!({\textbf r},{\textbf s},{\textbf t})\!))}{[\![R,S,T]\!]}.
$$

It is shown in \cite{IMW} that there exists the following commutative diagram
$$
\xyma{\frac{R\cap S\cap T}{[\![R,S,T]\!]} \ar@{=}[d] \ar@{->}[r] & 
\frac{{\textbf r}\cap {\textbf s}\cap {\textbf t}}{(\!({\textbf r},{\textbf s},{\textbf t})\!)}\ar@{->}[d] \\ \pi_2(\Omega X)\ar@{->}[r]& H_2(\Omega X)}
$$
Here $X$ is the homotopy pushout described above, $\Omega X$ is the loop space, the lower horizontal map is the Hurewicz homomorphism, the upper one is the natural map induced by $g\mapsto g-1$. For a connected space $Y$, the kernel of the second Hurewicz homomorphism $\pi_2(Y)\to H_2(Y)$ is a 2-torsion group \cite{IMW}. As a consequence we get the following theorem from \cite{IMW}. 

\begin{theorem}
For any group $G$ and its normal subgroups $R,S,T,$ the quotient
$$
\frac{G\cap (1+(\!({\textbf r},{\textbf s},{\textbf t})\!))}{[\![R,S,T]\!]}.
$$
is an abelian 2-torsion group. 
\end{theorem}

The author does not know how to prove this results without using the homotopy theory. 

What about the higher homotopy groups? Two normal subgroups give a possibility to model $\pi_2$ for a certain space, three normal subgroups correspond to $\pi_3$. Is it true that $n$ normal subgroups of a group allow to construct a space with a group-theoretical description of its $\pi_n$? The answer is "yes", however, under certain conditions. For a group $G$ and its normal subgroups $R_1,\dots, R_n,\ n\geq 2,$ construct the homotopy pushout of the $n$-dimensional cubical diagram with $2^n-1$ classifying spaces $K(G/\prod_{i\in I}R_i),\ I\subset \{1,\dots,n\}.$ There exist so-called connectivity conditions on the collection of subgroups $R_i,$ which imply that (see \cite{EllisMikhailov})
$$
\pi_n(\textrm{homotopy\ pushout})\simeq\frac{R_1\cap \dots\cap R_n}{[\![R_1,\dots, R_n]\!]}, 
$$
where
$$
[\![R_1,\dots, R_n]\!]:=\prod_{I\cup J=\{1,\dots, n\},\ I\cap J=\emptyset}[\bigcap_{i\in I}R_i,\bigcap_{j\in J}R_j]. 
$$

The main example we will consider here is the Wu formula for homotopy groups of $S^2.$ Let $F=F(x_1,\dots,x_n)$ be a free group of rank $n\geq 2$. Consider the following normal subgroups of $F$: $R_i=\langle x_i\rangle^F,\ i=1,\dots,n,\ R_{n+1}=\langle x_1x_2\dots x_n\rangle^F.$ The homotopy pushout of a diagram of corresponding $2^n-1$ classifying spaces is $S^2$. Therefore we get the following
$$
\pi_{n+1}(S^2)\simeq \frac{R_1\cap \dots\cap R_{n+1}}{[\![R_1,\dots, R_{n+1}]\!]}. 
$$
This is a version of the Wu formula, proved in \cite{Wu} using simplicial methods. In this particular case, the above commutator subgroup equals to the symmetric commutator subgroup
$$
[R_1,\dots, R_{n+1}]_S:=\prod_{\sigma\in \Sigma_{n+1}}[\dots [R_{\sigma(1)},R_{\sigma(2)}],\dots, R_{\sigma(n+1)}].
$$
Here $\Sigma_{n+1}$ is the group of $n+1$-permutations. That is, the homotopy groups of $S^2$ can be presented as 
$$
\pi_{n+1}(S^2)\simeq \frac{R_1\cap \dots\cap R_{n+1}}{[R_1,\dots, R_{n+1}]_S}.
$$
It turns out that the above intersection modulo the symmetric commutator coincides with the center of the quotient of the free group $F$ modulo the  
symmetric commutator, that is
$$
\pi_{n+1}(S^2)\simeq Z(F/[R_1,\dots,R_{n+1}]_S).  
$$
A generalization of this construction to the higher spheres as well as Moore spaces is given in \cite{MW}. For any $n, k>3,$ a finitely generated group $G_{n,k}$ given by explicit generators and relators is constructed, such that $\pi_n(S^k)\simeq Z(G_{n,k}).$ The group $G_{n,k}$ is defined in \cite{MW} as a certain quotient of the amalgamated square of the pure braid group on $n$ strands. 

Should we mention that the presentation of homotopy groups in this Section follows the idea of our formula $(\Phi)$? It is always an intersection of some subgroups modulo an obvious part. In this way, one can reflect on the difference between explicity of terms like $RS\cap RT$ and $R(S\cap T).$ 

\section{Wu-type formulas}
What are the first questions that come in mind when you look at Wu formula for homotopy groups of $S^2$? How to get some results on $\pi_*(S^2)$ using group theory? Can one find any new element from $\pi_*(S^2)$ using its group theoretic presentation? How to present generators of the known elements of the homotopy groups in terms of the free group? Here we will briefly discuss the last one. 

We know that $\pi_3(S^2)=\mathbb Z, \pi_4(S^2)=\mathbb Z/2, \pi_5(S^2)=\mathbb Z/2,\pi_6(S^2)=\mathbb Z/12,\dots$ all homotopy groups in degree $>1$ of $S^2$ are non-zero (see \cite{IMW1}). In general, the sequence of finite abelian groups $\pi_n(S^2), n\geq 4,$ is one of the most mysterious objects in math, it is difficult to speculate how far we are from its understanding. It is a strange luck that we can realize this extremely complicated sequence as a series of simply formulated subquotients of free groups. 

Looking at the Wu formula, one can also ask the following. What about other algebraic systems, different from groups? Clearly the same type of quotients can be considered for associative rings, Lie rings etc. As a rule, the answers will be easier. In the associative case this is almost obvious. The case of Lie rings (over $\mathbb Z$) is interesting and meaningful. Let $L_n=L(y_1,\dots, y_n),\ n\geq 2$ be a free Lie ring over $\mathbb Z.$ Consider its ideals $I_i=(y_i)^L,\ i=1,\dots, n,\ I_{n+1}=(y_1+\dots+y_n)^L.$ Define the Lie analog of the symmetric commutator of ideas
$$
[I_1,\dots, I_{n+1}]_S:=\prod_{\sigma\in \Sigma_{n+1}}[\dots [I_{\sigma(1)},I_{\sigma(2)}],\dots, I_{\sigma(n+1)}].
$$
There is the following isomorphism
\begin{equation}\label{e1}
\frac{I_1\cap \dots \cap I_{n+1}}{[I_1,\dots,I_{n+1}]_S}\simeq \bigoplus_{i\geq 1}E_{i,n}^1,
\end{equation}
where $E_{*,*}^1$ is the first page of the Curtis spectral sequence, defined via derived functors as
$$
E_{i,j}^1=L_j{\mathcal L}^i(\mathbb Z,1). 
$$
Here ${\mathcal L}^i$ is the $i$th Lie functor (see \cite{BM:2011} for the discussion of this spectral sequence and properties of derived functors). The values of $E_{*,*}^1$ are known and can be described in terms of Lambda-algebra (see \cite{6a}, \cite{Leib}). 

Now let's return to the problem of describing the generators in terms of free groups or Lie rings. It is clear that the generators can be chosen in different ways, since we work modulo the symmetric commutators. However we try to find ones with a simple form. Here are the results in low dimensions (see \cite{BartholdiMikhailov}). 

\vspace{.5cm}\noindent $\bf n=2.$ In this case, the generators are given as commutators $[x_1,x_2]$ in the group as well as $[y_1,y_2]$ the Lie ring case (here we will write the group-expressions in terms of $x$-s and Lie ring expressions in terms of $y$-s). Indeed, $[x_1,x_2]\in R_1\cap R_2\cap R_3\setminus [R_1,R_2,R_3]_S.$ In this case, $R_1\cap R_2\cap R_3=\gamma_2(F),\ [R_1,R_2,R_3]=\gamma_3(F)$ and $\gamma_2(F)/\gamma_3(F)\simeq \mathbb Z\simeq \pi_3(S^2).$

\vspace{.5cm}\noindent $\bf n=3.$ In this case, $\pi_4(S^2)=\mathbb Z/2,$ $E_{4,3}^1=\mathbb Z/2,\ E_{i,3}^1=0,\ i\neq 4.$ The generators of these $\mathbb Z/2$-terms are given as 
$$
[[x_1,x_2],[x_1,x_2x_3]],\ \ \ \ \ \ \ \ [[y_1,y_2],[y_1,y_3]].
$$
It can be easily checked that these terms lie in $R_1\cap\dots\cap R_4$ and $I_1\cap\dots\cap I_4$ respectively. Since we work in low dimensions it can be proved directly that they do not belong to the symmetric commutators. 

\vspace{.5cm}\noindent $\bf n=4.$ This case is just a suspension of the previous one, $\pi_5(S^2)=\mathbb Z/2,\ E_{8,4}^1=\mathbb Z/2,\ E_{i,4}^1=0,\ i\geq 8.$ The generators are given as
$$
[[[x_1,x_2],[x_1,x_2x_3]],[[x_1,x_2],[x_1,x_2x_3x_4]]]
$$
and 
$$
[[[y_1,y_2],[y_1,y_3]],[[y_1,y_2]],[y_1,y_4]]].
$$
The fact that these elements lie in the intersection of subgroups $R_i$ or ideals $I_i$ is obvious. To prove that they do not lie in the symmetric commutators we need some homotopy theory. It can be done using simplicial methods, by realizing these elements as cycles in the Milnor construction $F[S^1]$ and its Lie analog.  

\vspace{.5cm}\noindent $\bf n=5.$ This case is already complicated, $\pi_6(s^2)=\mathbb Z/12,$ $E_{6,5}^1=\mathbb Z/3, E_{8,5}^1=E_{16,5}^1=\mathbb Z/2,$ $E_{i,5}^1=0,\ i\neq 6,8,16.$ This horizontal line is the first place in the spectral sequence where one can see a non-trivial gluing of $E^\infty$-term: two cells in degrees 8 and 16 with values $\mathbb Z/2$ are glued into $\mathbb Z/4\subset \pi_6(S^2).$ It is easy to write down the generator of $E_{16,5}^1,$
since it comes as a suspension of the previously written element $E_{8,4}^1$, namely 
$$
[[[[y_1,y_2],[y_1,y_3]],[[y_1,y_2]],[y_1,y_4]]],[[[y_1,y_2],[y_1,y_3]],[[y_1,y_2]],[y_1,y_5]]]].
$$
The term $E_{6,5}^1=\mathbb Z/3$ is generated by the element (see \cite{BartholdiMikhailov} for the proof using simplicial methods and derived functors)
\begin{align*}
& \alpha_3:=[[[y_1, y_5], [y_2, y_5]], [y_3, y_4]] - [[[y_1, y_5], [y_3, y_5]], [y_2, y_4]]\\ 
& + [[[y_1, y_5], [y_4, y_5]], [y_2, y_3]] + [[[y_2, y_5], [y_3, y_5]], [y_1, y_4]]\\
& - [[[y_2, y_5], [y_4, y_5]], [y_1, y_3]] + [[[y_3, y_5], [y_4, y_5]], [y_1, y_2]]
\end{align*}
The term $E_{8,5}^1$ as well as group-case liftings are much more complicated (see \cite{BartholdiMikhailov}). For example, the 3-torsion from $\pi_6(S^2)$ can be written as a product of 14 commutators in a free group of weight $\geq 6.$

The element $\alpha_3$ corresponds to the Serre element $\mathbb Z/3\subset \pi_6(S^2)$. Analogous picture takes a place for all primes. For an odd prime $p$, the Serre element of order $p$ appear in the homotopy group $\pi_{2p}(S^2)$. These elements can be easily seen from the structure of the first page of the spectral sequence $E_{2p,2p-1}^1=\mathbb Z/p$ [these terms are labeled as $\lambda_1$ in the language of Lambda-algebras]. It turns out that the $E_{2p,2p-1}^1$-term is isolated from other $\mathbb Z/p$-torsion terms of the spectral sequence, hence $E_{2p,2p-1}^1=E_{2p,2p-1}^\infty$ [In the theory of spectral sequences such kind of arguments sometimes are called \textit{lacunary reasons}]. Let $y_i$ for $i=1,\dots,2p-1$ be free generators of a free Lie
  algebra, consider the following element
  \begin{multline*}
    \alpha_p=\sum_{\kern-1cm\substack{\rho\in\Sigma_{2p-2}\text{ a $2^{p-1}$-shuffle}\\\rho(1)<\rho(3)<\dots<\rho(2p-5)}\kern-1cm}(-1)^\rho[[y_{\rho(0)},y_{2p-2}],[y_{\rho(1)},y_{2p-2}],[y_{\rho(2)},y_{\rho(3)}],\dots, [y_{\rho(2p-4)},y_{\rho(2p-3)}]];
  \end{multline*}
  the sum is taken over all permutations
  $(\rho(0),\dots,\rho(2p-3))\in\Sigma_{2p-2}$ satisfying
  $\rho(0)<\rho(1),\dots,\rho(2p-4)<\rho(2p-3)$ as well as
  $\rho(1)<\rho(3)<\dots<\rho(2p-5)$. Here we use the left-normalized notation, i.e. $[x,y,z]:=[[x,y],z]$. Then $\alpha_p$ presents a generator of $L_{2p-1}\mathcal L^{2p}(\mathbb Z,1)=E_{2p,2p-1}^1$ (see \cite{BartholdiMikhailov} for the proof). 

Looking at the elements of free groups and Lie rings like $[[x_1,x_2],[x_1,x_2x_3]]$ or $\alpha_3$ one can get some impression of \textit{homotopy patterns}. Next we will see how they work in the context of dimension subgroups. 

\section{Classical dimension subgroups} Is it true that, for any group $G$, and $n\geq 1$, $D_n(G)=\gamma_n(G)$? This question is known as \textit{dimension problem} and has a long history. For a detailed discussion of this problem we refer to \cite{Passi}, \cite{Gupta:1987}, \cite{MP:2009}. The first results in this direction are due to Magnus and Witt. They proved that, for a free group $F$, the dimension subgroups coincide with lower central series \cite{Magnus:1937}, \cite{Witt:1937}. Incorrect proofs of the dimension problem appeared more than once, see
\cite{Cohn}, \cite{Losey} and [\cite{MKS}, {Theorem~5.15(i)}]. The first example of a group with $D_4(G)\neq \gamma_4(G)$ is due to Ilya Rips \cite{Rips}. The group constructed in \cite{Rips} has the order $2^{38}$ and it seems that this is the smallest finite group with the property $D_4\neq \gamma_4.$ The next point we have to mention as a the history of the question is the series of works of Narain Gupta. In order to describe the dimension subgroups and solve the dimension problem completely, Gupta spent about 20 years developing a special calculus. As a final result he published the paper \cite{Gupta:2002}, where he claims that the dimension property holds for all groups of odd order. In particular, it follows from his claim that, for an odd prime $p$, it is not possible to construct a group $G$ and $n\geq 1$ with $D_n(G)/\gamma_n(G)\supseteq \mathbb Z/p.$ In fact, Gupta claimed even more, that, for any group $G$, the dimension quotients $D_*(G)/\gamma_*(G)$ are just $\mathbb Z/2$-vector spaces. The last statement was written in the non-published manuscript of Gupta, which was available from 90-s to the experts in the area.  

The proofs given in the mentioned papers of Gupta are extremely complicated. During many years the author together with I.B.S. Passi tried to understand these proofs. It was clear already about 10 years ago that they contain gaps, however, it was not easy to find counterexamples to the main statements. Finally, the following result is proved in \cite{BartholdiMikhailov}: 

\begin{theorem} For any prime $p$, there exists a group $G$ and integer $n$, such that $D_n(G)/\gamma_n(G)$ contains $\mathbb Z/p$ as a subgroup. 
\end{theorem}

Among other things, a small finite group $G$ with $D_7(G)\neq \gamma_7(G)$ is constructed in \cite{BartholdiMikhailov}. The needed statement that $D_7(G)\neq \gamma_7(G)$ is checked using GAP. The order of $G$ is $3^{494}$ and this is a 3-group without dimension property of the smallest order the authors were able to construct. 

In \cite{BartholdiMikhailov}, both categories are considered, groups and Lie rings. What often happens, the computations for Lie rings are simpler. Here we will give detailed examples for the case of Lie rings. For a Lie ring over integers $L$, consider its universal enveloping algebra $U(L)$. The algebra $U(L)$ admits the augmentation ideal $\omega.$ The Lie ring $L$ embeds in $U(L),$ and the dimension subgroups of $L$ are defined as $\delta_n(L):=L\cap \omega^n.$ The lower central series term $\gamma_n(L)$ lies in $\delta_n(L)$ and almost all main statements of the theory of dimension subgroups can be extended from group to Lie rings, with simpler proofs (see \cite{BP}). In particular, for any $L$, $\delta_n(L)=\gamma_n(L), n=1,2,3$ and there exists a Lie analog of Rips example such that $\delta_4/\gamma_4=\mathbb Z/2.$ The following result also is from \cite{BartholdiMikhailov}

\begin{theorem}For any prime $p$, there exists a Lie ring $A$ and integer $n$, such that the abelian group $\delta_n(A)/\gamma_n(A)$ contains $\mathbb Z/p$ as a subgroup. 
\end{theorem}

In a Lie ring presentation, we introduce the
following notation: for $d\in\mathbb N$, when we write a generator $y^{(d)}$
of degree $d$ we mean a list of generators $y_1,\dots,y_d$; when $y^{(d)}$
the left-normed iterated
commutator $y^{(d)}:=[y_1,\dots,y_d]$. Thus, for example, `$\langle y_1^{(2)},y_2^{(3)}\mid [y_1,y_2]\rangle$'
is shorthand for
`$\langle
y_{1,1},y_{1,2},y_{2,1},y_{2,2},y_{2,3}\mid[[y_{1,1},y_{1,2}],[y_{2,1},y_{2,2},y_{2,3}]]\rangle$'. The following is proved in \cite{BartholdiMikhailov}. 
 Given an integer $s\ge3$, there are integers $e,c_0,\dots,c_s$ and
  $n=c_0+\dots+c_s$ such that, for the Lie ring
  \[A=\langle y_0\dots,y_s,z_0^{(c_0)},\dots,z_s^{(c_s)}\mid y_0+\dots+y_s=0, e^{c_i}y_i=z_i\text{ for }i=0,\dots,s\rangle
  \]
  there exists a natural embedding
  \[\bigoplus_i E_{i,s}^1\hookrightarrow\delta_n(A)/\gamma_n(A).\]
  
To illustrate how does it work, we first rewrite the formula (\ref{e1}) as follows. Take $L$ a free Lie ring with generators $y_0,\dots, y_s$ and one relation $y_0+\dots+y_s=0.$ Set $I_i=(y_i)^L.$ In this notation, 
\[
\frac{I_0\cap \dots \cap I_s}{[I_0,\dots,I_s]_S}\simeq \bigoplus_{i\geq 1}E_{i,s}^1.
\]
  
Consider the universal enveloping algebra $U(L)$, the corresponding ideals $J_i=y_i U(L)$ in $U(L)$, and their symmetric
product:
\[
(J_0,\dots, J_s)_S=\sum_{\rho\in\Sigma_{s+1}}J_{\rho(0)}\cdots J_{\rho(s)}.
\]
The natural map $L\to U(L)$ induces
  \[\xyma{
      \frac{I_0\cap\dots\cap I_s}{[I_0,\dots,I_s]_S} \ar@{=}[d]\ar@{->}[r] & \frac{J_0\cap\dots\cap J_s}{(J_0,\dots,J_s)_S} \ar@{=}[d]\\
      \bigoplus_{i\ge1} E_{i,s}^1\ar@{->}[r] & H_s(U(L[S^1])).}
  \]
Here $U(L[S^1])$ is the universal enveloping algebra of the simplicial Lie ring $L[S^1],$ it has infinite cyclic homology groups in all dimensions. At the same time, $E_{i,j}^1$-terms of the lower central
  series spectral sequence for $S^2$ are finite for all $j\geq 3$. It follows that
  the map is $0$. Therefore, for $s\geq 3$ we have
  $I_0\cap\dots\cap I_s\le L\cap(J_0,\dots, J_s)_S$ when
  considered in the universal enveloping algebra. See \cite{BartholdiMikhailov} for details. 

Next we add the relations $e^{c_i}y_i=z_i$, for a fixed choose of $e, c_i, i=0,\dots,s.$ Remind that the elements $z_i$ are of degrees $c_i$. These relations imply that, for any $w\in I_0\cap\dots\cap I_s$, the image of $e^nw$ in $A$, lies in $\delta_n.$ An element $w$ is our \textit{homotopy pattern}, we see that, if we write it in the universal enveloping algebra $U(L),$ it lies in the symmetric product $(J_0,\dots, J_n)_S,$ at the same time, it does not lie in the symmetric product $[I_0,\dots, I_s]_S.$ Now the relations $e^{c_i}y_i=z_i$ guarantee that $e^nw$ lies in the $n$th augmentation power of the universal enveloping algebra. There is no obvious reason for the element $e^nw$ to lie in the $n$th term of lower central series of $A$. The detailed analysis shows that one can choose the constants $e, c_i,$ such that $e^nw$ will not lie in $\gamma_n(A)$ (see \cite{BartholdiMikhailov}). We can call the described method by \textit{bloating of a homotopy pattern}. 

In particular, we can take $e=p,$ $s=2p-1$ and realize the $p$-torsion elements $\alpha_p$ described in the previous section inside dimension quotients of Lie rings. 

 In some cases it is possible to find $c_i$-s sufficiently small. The following two examples are from \cite{BartholdiMikhailov}, they are checked using computer assistance. Consider the Lie ring
\begin{multline*}
  A=\langle y_0,y_1,y_2,y_3,z_0^{(1)},z_1^{(2)},z_2^{(2)},z_3^{(2)}\mid\\
  y_0+y_1+y_2+y_3=0,\;y_0=2^6z_0,\;2^6y_1=2^5z_1,\;2^5y_2=2^3z_2,\;2^3y_3=z_3\rangle
\end{multline*}
and the element $\omega=[[y_0,y_1],[y_0,y_2]]$. In that Lie algebra,
we have $\omega\in\delta_{7}(A)\setminus\gamma_{7}(A)$ and
$2\omega\in\gamma_{7}(A)$. So, our \textit{homotopy pattern} which generated $L_3\mathcal L^4(\mathbb Z,1)$ and $\pi_4(S^2),$ makes the difference between $\delta_7$ and $\gamma_7.$

For $p=3$ the situation is similar and the construction can be simplified. Consider the following Lie ring
\begin{align*}
A=\langle y_{ij},z_{ij}^{(i+j+1)}\text{ for }0\le i<j\le 5\ |\ 
  & y_{01}+y_{02}+y_{03}+y_{04}+y_{05}=0,\\
  & {-}y_{01}+y_{12}+y_{13}+y_{14}+y_{15}=0,\\
  & {-}y_{02}-y_{12}+y_{23}+y_{24}+y_{25}=0,\\
  & {-}y_{03}-y_{13}-y_{23}+y_{34}+y_{35}=0,\\
  & {-}y_{04}-y_{14}-y_{24}-y_{34}+y_{45}=0,\\
  & 3^{i+j}y_{ij}=z_{ij}\text{ for }0\le i<j\le 5\rangle
\end{align*}
Then the element $\omega=3^{15}([y_{04},y_{14},y_{23}]-[y_{04},y_{24},y_{13}]+[y_{04},y_{34},y_{12}]+[y_{14},y_{24},y_{03}]-[y_{14},y_{34},y_{02}]+[y_{24},y_{34},y_{01}])$ belongs to $\delta_{18}(A)\setminus\gamma_{18}(A)$ with $3\omega\in \gamma_{18}(A).$ One can easily recognize our \textit{homotopy pattern} $\alpha_3$ rewritten in the way $[y_i,y_j]\rightsquigarrow y_{i-1,j-1}.$

The situation for groups is similar. We start with a free group \[F=\langle x_0,\dots, x_s\ |\ x_0\dots x_s=1\rangle,\ R_i=\langle x_i\rangle^F,\]
take an element $q$ from the intersection $R_0\cap \dots \cap R_s\setminus [R_0,\dots, R_s]_S.$ Since, for $s\geq 3,$ the Hurewicz homomorphism $\pi_{s+1}(S^2)=\pi_s(\Omega S^2)\to H_s(\Omega S^2)$ is the zero map, the element $q-1$, written in the group ring $\mathbb Z[F],$ belongs to the symmetric product of ideals $(r_0,\dots, r_s)_S.$ Next we use a \textit{bloating of a homotopy pattern} applied to the group-case. The difference between Lie ring and group cases is that, in the case of groups, we can not transfer exponents in a free way: $[x^e,y]\neq [x,y^e].$ That is, the process have to be done more carefully. Again we refer to \cite{BartholdiMikhailov} for details. The $\mathbb Z/p$-torsion terms in the dimension quotients can be realized, for example, using Serre elements $\mathbb Z/p\subset\pi_{2p}(S^2).$ The experiments with GAP show that Rips-type examples with $D_4\neq \gamma_4$ also can be constructed using the described method, bloating of a homotopy pattern $[[x_1,x_2],[x_1,x_2x_3]]$ which corresponds to the homotopy element $\mathbb Z/2=\pi_4(S^2)$.

\section{Limits. Speculative functor theory.}
\textit{"All Rajayoga, for instance, depends on this perception and experience that our inner elements, combinations, functions, forces,
can be separated or dissolved, can be new-combined and set to
novel and formerly impossible workings or can be transformed
and resolved into a new general synthesis by fixed internal
processes." (Sri Aurobindo, "The Synthesis of Yoga")}

A standard way to define a group or Lie algebra, is combinatorial, i.e. via generators and relations. In order to construct an algebraic objects with complicated properties, one can play with generators and relations. Take, for example, two symbols "a" and "b" and one relator 
\[
\langle a,b\ |\ a^{-1}b^2ab^{-1}a^3b^{-1}=1\rangle.
\]
It turns out that the resulting group has interesting and non-obvious properties. When we consider some \textit{functor}, we usually mean that this functor comes from a natural consideration and is not constructed in a speculative way. There are no obvious combinatorial games which give a possibility to define functors in terms of generators and relations. However, there is one non-obvious way, this section will be about it. We will see that the characters we discussed before, like derived functors, dimension quotients and group homology will appear in that functorial constructions. 

Let $G$ be a group. By $\textsf{Pres}(G)$ we denote the \textit{category of presentations} of $G$ with objects being free groups $F$ together with epimorphisms to $G$.
Morphisms are group homomorphisms over $G$. For a functor $\mathcal F:\textsf{Pres}(G)\to \textsf{Ab}$ from the category $\textsf{Pres}(G)$ to the category of abelian groups, one can
consider the (higher) limits $\ilimit^i\, \mathcal F,\ i\geq 0,$ over the category of presentations. That is, we fix our group $G$, consider free presentations $R\hookrightarrow F\twoheadrightarrow G$ and make functorial (on $F,R$) constructions $\mathcal F(F,R).$ The limits $\ilimit^i\, \mathcal F(F,R),\ i\geq 0,$ will depend only on $G$ and, moreover, present functors from the category of all groups to abelian groups. 

The category $\textsf{Pres}(G)$ is strongly connected and has pair-wise coproducts. The limit $\ilimit F=\ilimit^0 F$ has the following properties. For any $c\in \textsf{Pres}(G),$
\begin{equation}\label{equa}
\ilimit\mathcal F=\{x\in \mathcal F(c)|\forall\,c'\in\textit{Pres}(G), \ \varphi, \psi:c\to c',\ \mathcal F(\varphi)(x)=\mathcal F(\psi)(x)\}.
\end{equation}
Moreover, it is known \cite{MP2019} that the limit of a functor from a strongly connected category with pair-wise coproducts is equal to the equalizer
$$\ilimit \mathcal F \cong \textrm{eq}(\mathcal F(c)\rightrightarrows \mathcal F(c\sqcup c))$$
for any $c\in \textsf{Pres}(G)$.
 In particular, this equalizer does not depend on $c$.

The standard and one of the simplest examples is the following. For a group $G=F/R,$ the Hopf formula for the second homology is very useful: $H_2(G)=\frac{R\cap \gamma_2(F)}{[F,R]}.$ The second homology can be presented as limit as well:
\[
H_2(G)=\ilimit R/[F,R]. 
\]
That is, the quotient $\frac{R\cap \gamma_2(F)}{[F,R]}$ is the maximal subgroup of $R/[F,R],$ which depends on $G$ only. 

The limits $\ilimit^i\, \mathcal F$ are
studied in the series of papers \cite{EM}, \cite{IvanovMikhailov2015}, \cite{MP:2015a}, \cite{IMP}, \cite{MP2019}. Next we will give examples which illustrate the variety and complexity of functors which one can obtain playing with limits. 

The group homology can be presented as follows (see \cite{EM}, \cite{IvanovMikhailov2015})
\begin{align*}
& \ilimit (R_{ab}^{\otimes n})_F=H_{2n}(G), n\geq 1\\ 
& \ilimit^1 (R_{ab}^{\otimes n})_F=H_{2n-1}(G). 
\end{align*}
Here the tensor powers of the relation modules $R_{ab}$ are considered with diagonal action of $F.$ The derived functors $L_1\textsf S^2$ and $L_1\textsf S^3$ which we discussed in Section 2, can be described via limits (see \cite{MP:2015a}, \cite{MP:2016})
\begin{align*}
& \ilimit \frac{\gamma_2(F)}{\gamma_2(R)\gamma_3(F)}=L_1\textsf S^2(G_{ab}),\\     
& \ilimit \frac{\gamma_3(F)}{[\gamma_2(R),F]\gamma_4(F)}=L_1\textsf S^3(G_{ab}),\\
& \ilimit \frac{\gamma_3(F)}{\gamma_3(R)\gamma_4(F)}=L_2 {\textsf L}_s^3(G_{ab}),\\ 
& \ilimit \frac{\gamma_4(F)}{[\gamma_2(R),F,F]\gamma_2(\gamma_2(F))\gamma_5(F)}=L_1\textsf S^4(G_{ab}),\\ 
& \ilimit \frac{\gamma_2(R)}{\gamma_2([R,F])\gamma_3(R)}=L_1\textsf S^2(H_2(G)). 
\end{align*}
Comparing the first three limits with the results of Section 2, we see that, in certain cases, when the formula $(\Phi)$ can be applied, our \textit{homotopy patterns} can be described as limits
\[
\ilimit \frac{\textsf{Whole structure}}{\textsf{Obvious part of an intersection}}=\textsf{Homotopy pattern}.
\]
Observe that, a simple deformation of the considered functor $\textsf{Pres}\to \textsf{Ab}$ may change the limits completely. For example, 
\[
\ilimit \frac{\gamma_3(F)}{[\gamma_2(F),R]\gamma_4(F)}=0. 
\]
Let's mention a couple of exotic examples. 
If a group $G$ does not have a 2-torsion, then \cite{IvanovMikhailov2015}
\begin{align*}
& \ilimit \gamma_2(R)/[\gamma_2(R),F]=H_4(G; \mathbb Z/2),\\ 
& \ilimit^1 \gamma_2(R)/[\gamma_2(R),F]=H_3(G; \mathbb Z/2). 
\end{align*}
Recall the Fox quotient $\frac{F(3,R)}{G(3,R)}$ from Section 2. This quotient depends on $F$ and $R,$ not only on $G$. The limit of this quotient is computed in \cite{MP2019}:
\[
\ilimit \frac{F(3,R)}{G(3,R)}=L_1\textsf S^2(L_1\textsf S^2(G_{ab})). 
\]

Next result shows how to present the fourth dimension quotient via limits (see \cite{MP2019} for the proof)
\begin{theorem}\label{th61}
There is a natural short exact sequence
\[
\ilimit \frac{R\cap \gamma_2(F)}{\gamma_2(R)(R\cap \gamma_4(F))}\hookrightarrow \ilimit \frac{\gamma_2(F)}{\gamma_2(R)\gamma_4(F)}\twoheadrightarrow \frac{D_4(G)}{\gamma_4(G)}.
\]
\end{theorem}
\vspace{.5cm}
 Theorem \ref{th61} shows how to describe the fourth dimension quotient as a functor without using the group ring. The limits in Theorem \ref{th61} are just equalizers (\ref{equa}), one can compute them for simple examples. The author thanks L. Bartholdi for computing these limits for Rips-type examples using computer assistance. These computations show that the short exact sequence from Theorem \ref{th61} does not split.  

We finish the paper by briefly reviewing so-called $\textbf{fr}$-\textit{language} (see \cite{IM2017}, \cite{IMP}). 
The ideals $\textbf f=(F-1)\mathbb Z[F],\ \textbf r=(R-1)\mathbb Z[F]$ define functors $\textsf{Pres}(G)\to \textsf{Ab}.$ Moreover, all possible product of ideals $\textbf{f}, \textbf{r}$, their sums and intersections define fucntors  $\textsf{Pres}(G)\to \textsf{Ab}$ as well, and we can ask how to describe their limits. One can take any \textit{sentence} on symbols $\textbf{f}, \textbf{r}$ like $\textbf{rr}+\textbf{ffr}+\textbf{frf}+\textbf{rff}$ or $\textbf{rrrf}+\textbf{frrr}$ and consider their $\ilimit^i$ as functors. The author does not know any unified method of description of limits for a given $\textbf{fr}$-sentence. For any particular case there are some special tricks, based on homological algebra or group theory. Sometimes the results are surprising. For example, $\ilimit^2$ of two mentioned sequences are the following:  
\begin{align*}
& \ilimit^2 (\textbf{rr}+\textbf{ffr}+\textbf{frf}+\textbf{rff})=G_{ab}\otimes G_{ab},\\ 
& \ilimit^2 (\textbf{rrrf}+\textbf{frrr})=H_5(G). 
\end{align*}
Here are some more examples of computations, which show that the variety of functors which can be presented as limits of $\textbf{fr}$-sentences is rich enough: 
\begin{align*}
& \ilimit^1(\textbf{rff+frr})=\textsf{Tor}(H_2(G),G_{ab}),\\
& \ilimit^1(\textbf{rr+frf+rff})=H_2(G,G_{ab}),\\
& \ilimit^1(\textbf{rr+frf})=H_3(G),\\
& \ilimit^2(\textbf{rr+frf})={\textbf g}\otimes_{\mathbb Z[G]}{\textbf g},\\ 
& \ilimit^1(\textbf{rrf+frr})=H_4(G),\\ 
& \ilimit^1(\textbf{rr+fff})=\textsf{Tor}(G_{ab},G_{ab}),\\ 
& \ilimit^2(\textbf{rr+fff})=G_{ab}\otimes G_{ab}.
\end{align*}

The point of this theory (which we also call $\textbf{fr}$-\textit{language}), is that the formal manipulations with codes in two letters may induce deep and
unexpected transformations of functors. Simple transformations of $\textbf{fr}$-codes, like changing the symbol $\textbf r$ by $\textbf f$ in a certain place, adding a monomial to the $\textbf{fr}$-code etc,
induce natural transformations
of (higher) limits determined by these $\textbf{fr}$-codes. For example, the transformation of the $\textbf{fr}$-codes
$$
\textbf{rr+frf}\rightsquigarrow \textbf{rr+frf+rff}
$$
induces the natural transformation of functors
$$
H_3(G)=\ilimit^1(\textbf{rr+frf})\rightsquigarrow \ilimit^1(\textbf{rr+frf+frr})=H_2(G,G_{ab}).
$$
Here the map $H_3(G)\to H_2(G, G_{ab})$ is constructed as \[ H_3(G)=H_2(G,{\textbf g})\to H_2(G,{\textbf g}/{\textbf g}^2)=H_2(G,G_{ab}),\] where the last map is induced by the natural projection
${\textbf g}\twoheadrightarrow {\textbf g}/{\textbf g}^2=G_{ab}.$

We end this Section with an observation that, in many cases, when the formula $(\Phi)$ can be applied to the structures described in terms of $\textbf{f}$ and $\textbf{r},$ the \textit{homotopy patterns} can be seen via limits. For example, the well-known description of the $2n$-th homology $(n\geq 1)$, \[H_{2n}(G)=\frac{\textbf r^n\cap \textbf{fr}^{n-1}\textbf f}{\textbf r^n\textbf f+\textbf f\textbf{r}^n}, \]
represents the formula $(\Phi).$ Simple computation (see \cite{IM2017}) shows that 
\[
\ilimit^1 (\textbf r^n\textbf f+\textbf f\textbf{r}^n)=H_{2n}(G), n>1
\]
The case $n=1$ is an exception: $\ilimit^1 (\textbf{rf}+\textbf{fr})={\textbf g}\otimes_{\mathbb Z[G]}{\textbf g}\supset H_2(G).$ 
In such kind of cases, the higher limits give a way to consider derived versions of the \textit{homotopy patterns} as well. It seems that this is a good point to end this survey.

\end{document}